# The Blackwell-Cover Guessing Strategy for Arbitrary Probability Spaces

By James D. Stein

Department of Mathematics and Statistics, California State University (Long Beach)

**Abstract:** The Blackwell-Cover Guessing Strategy is well known in the case where the guesser is confronted with two envelopes containing different sums of money. Guesser is allowed to open one and count the money in it and can then use a random variable to guess with success probability greater than ½ which envelope contains the larger amount. It is less well known in the case where the guesser is trying to guess which of two coins has the greater heads probability. Guesser flips one, if it land heads guess that it has the larger heads probability, if it lands tails guess that the other coin has the larger heads probability. Both of these entail a single observation. We show that there is a similar Blackwell-Cover Guessing Strategy for arbitrary probability spaces, which enables guesser to guess which of two subsets in a σ-field has the larger probability based on the observation of a random point in the sample space, thus potentially greatly expanding the range of applications of the Blackwell-Cover Guessing Strategy. We extend this result to the case where we are trying to guess which of N subsets in a σ-field has the largest probability based on the observation of a random point in the sample space; we can guess correctly with probability > 1/N.

We also examine two different well-known two-observation strategies in the case where we are trying to guess which of two subsets in the σ-field has the larger probability. One yields exactly the same probability of a successful guess as the single-observation guessing strategy. The other neither dominates nor is dominated by the single observation guessing strategy. As a result, there is no inherent advantage in making two observations rather than one.

**Section I -The Blackwell-Cover Guessing Strategy with a Single Observation**

The only surprising thing about the following theorem is that nobody seems to have noticed it.

**Theorem 1** – Let S be a set, and let $\mathscr{F}$ be a σ-field of subsets of S, P a probability measure on $\mathscr{F}$. If we choose x randomly in S and E $\epsilon$ $\mathscr{F}$, P(E) denotes the probability that x will belong to E. Given two sets A and B in $\mathscr{F}$ with $P(A) \neq P(B)$, define the Blackwell-Cover Guessing Strategy as follows. Choose one of the two sets A and B equiprobably, and choose x randomly. If x belongs to the chosen set, guess that that set has the larger probability, otherwise guess that the other set has the larger probability. The guess will be successful with probability $\frac{1}{2} + \frac{1}{2} |P(A) - P(B)|$.

**Proof:** Assume that $P(A) > P(B)$. With probability ½, the equiprobably chosen set will be A. The guessing strategy dictates that we guess correctly if x $\epsilon$ A; the combined probability that we guess correctly is therefore ½ P(A).

With probability ½, the equiprobably chosen set will be B. The guessing strategy dictates that we guess correctly if x $\epsilon$ B', the complement of B. The combined probability that we guess correctly is therefore ½ (1-P(B)).

The probability of a correct guess is therefore $\frac{1}{2} P(A) + \frac{1}{2} (1 - P(B)) = \frac{1}{2} + \frac{1}{2} (P(A) - P(B))$. If $P(B) > P(A)$, the probability of a correct guess would be $\frac{1}{2} + \frac{1}{2} (P(B) - P(A))$. The theorem follows. █

**Cor. 1.1 -** Let S be a set and let $\mathscr{F}$ be a σ-field of subsets of S, P a probability measure on $\mathscr{F}$. Let $A_1, \ldots, A_N$ belong to $\mathscr{F}$ and have different probabilities, we describe a Blackwell-Cover Guessing Strategy which guesses the set $A_k$ having maximum value of the $P(A_j)$ for $1 \leq j \leq N$, and guesses successfully with probability > 1/N. Choose one of the sets $A_j$ equiprobably, and choose x randomly in S. If x belongs to the chosen set, guess that it has the maximum value of $P(A_j)$, otherwise guess one of the other sets randomly uniformly.. If $P(A_k)$ has the maximum value, the guess will be successful with probability $\frac{1}{N}(1 + P(A_k) - \text{average}\{P(A_j)|j \neq k\})$.

**Proof:** With probability 1/N, the set $A_k$ will be chosen equiprobably, and the combined probability that $A_k$ will be chosen equiprobably and the guess will be successful is $P(A_k)/N$. If $j \neq k$, the combined probability that $A_j$ will be chosen, x will not belong to $A_j$, and the random uniform choice of the other N-1 sets will select $A_k$ is $\frac{1}{N}(1 - P(A_j))\frac{1}{N-1}$. So the probability that the guess will be successful using this strategy is

$$\frac{1}{N}\left(P(A_k) + \frac{1}{N(N-1)}\sum_{j\neq k}(1 - P(A_j))\right) = \frac{1}{N}(1 + P(A_k) - \text{average}\{P(A_j)|j \neq k\})$$

This completes the proof. █

**Example 1 (Which is larger, Africa or South America?)** – Suppose one is asked to guess which continent has the largest surface area – Africa or South America. Use the area model for probability, assuming that the set S is the entire surface of the Earth. Pick one of the two continents (assume for the purpose of description we have chosen Africa) and get a globe of the Earth and a dart. Spin the globe and toss the dart at the globe. If it sticks in Africa, guess that Africa is the larger continent, if not, guess that South America is the larger continent.

**Example 2 (Are you more likely to encounter a left-handed magician or a blue-eyed juggler?)** – In this case, the sample space can be viewed as a collection of binary strings constituting yes/no responses to questions on a questionnaire, one string for each respondent. The σ-field is constructed by taking unions and intersections of each of the sets corresponding to a specific string, and the probability of that set is simply its cardinality divided by the total number of respondents. To apply the Blackwell-Cover Guessing Strategy, select one of the two sets (left-handed magicians or blue-eyed jugglers) equiprobably – assume for the purpose of description that we have selected left-handed magicians. Now select a random respondent. If that respondent is a left-handed magician, guess that left-handed magicians are more prevalent than blue-eyed jugglers, and conversely.

**Section II – Blackwell-Cover Guessing Strategies for Two Sets Given Two Observations**

Any single-observation Blackwell-Cover Guessing Strategy can be described by the phrase "stay on success, switch on failure." What is "staying" or "switching" is the initial equiprobable choice of which is larger. This is an envelope (as in the original Blackwell-Cover Guessing Strategy), a coin (the Blackwell-Cover Guessing Strategy for Bernoulli trials), or a set in a σ-field (this paper).

In this section we examine two different guessing strategies using two observations (picking two elements from S), continuing the "stay or switch" theme. Throughout this section, we let a = P(A) and b = P(B). We assume WLOG that a > b.

**Strategy 1 – Stay on Success**

For the first strategy, which we abbreviate by the phrase Stay on Success, we first equiprobably select one of the two sets, and make two observations. If both observations yield elements belong to that set, we guess that it is the set with the larger probability. If both elements lie outside that set, we guess that the other set is the one with the larger probability. If one element lies inside the set and the other lies outside the set, we flip a fair coin to guess which set has the larger probability.

**Theorem 2** – The success probability of Stay on Success is ½ (1 + a – b).

**Proof:** We choose the set A with probability ½. The probability that both elements will belong to A is $a^2$, in which case the guess is correct. The probability that precisely one of the two elements will belong to A is 2a(1-a), so the probability that we will guess correctly in this case is a(1-a). The combined probability of choosing A and guessing correctly is therefore ½ ($a^2$ +a(1-a)) = ½ a.

We now choose the set B with probability ½. The probability that both elements will not belong to B is $(1-b)^2$, in which case the guess is correct. The probability that precisely one of the two elements will belong to B is 2b(1-b), so the probability that we will guess correctly in this case is b(1-b). The combined probability of choosing B and guessing correctly is therefore seen to be ½ ($(1-b)^2$ +b(1-b)) = ½ (1-b).

The combined probability of using this strategy is therefore ½ a + ½ (1-b). ∎

It is somewhat surprising that this strategy has exactly the same success probability as the single-observation Blackwell-Cover Guessing Strategy described in the previous section.

**Strategy 2 - Switch on Failure**

In this strategy, we choose one of the two sets equiprobably and perform the first of the two observations. If the first element observed belongs to the chosen set, we then randomly observe the second element. If the second element also belongs to the chosen set, we guess that the chosen set has the larger probability. If the second element does not belong to the chosen set, we flip a fair coin to guess which set has the larger probability.

If the first element does not belong to the chosen set, we "switch on failure", switching our choice of set to the other set. If the second element belongs to the other set, we guess that that set is the one with the higher probability. If the second element does not belong to the other set, we flip a fair coin to guess which set has the larger probability.

**Theorem 3** – The success probability of Switch on Failure is $\frac{1}{2} + \frac{1}{4}(a - b)(1 + a + b)$.

**Proof:** In the table below, 1st choice refers to the initial equiprobable choice of set. 2nd choice is the same set if the first element observed belongs to it, and is the other set if the first element observed does not belong to it. **Comb Prob** is the combined probability of making the first choice of set and having the first and second elements observed in the locations described in the 2nd and 4th columns in the table below. **Prob of Correct Guess** is the probability of making the correct guess as to the set with the larger probability, given the conditions described in that row.

| 1st choice | 1st element | 2nd choice | 2nd element | **Comb Prob** | **Prob of Correct Guess** |
|---|---|---|---|---|---|
| A | in A | A | in A | $\frac{1}{2} a^2$ | 1 |
| A | in A | A | not in A | $\frac{1}{2} a(1-a)$ | ½ |
| A | not in A | B | in B | $\frac{1}{2} (1-a)b$ | 0 |
| A | not in A | B | not in B | $\frac{1}{2} (1-a)(1-b)$ | ½ |
| B | in B | B | in B | $\frac{1}{2} b^2$ | 0 |
| B | in B | B | not in B | $\frac{1}{2} b(1-b)$ | ½ |
| B | not in B | A | in A | $\frac{1}{2}(1-b)a$ | 1 |
| B | not in B | A | not in A | $\frac{1}{2} (1-b)(1-a)$ | ½ |

The probability of success for this strategy is therefore

$$\frac{1}{2} a^2 + \frac{1}{4} (a(1-a) +(1-a)(1-b) + b(1-b) +(1-b)(1-a)) + \frac{1}{2} (1-b)a$$

After some simplification, this reduces to ½ + ¼ (a - b)(1 + a + b). ■

The single observation success probability for the Blackwell-Cover Guessing Strategy is ½ + ½ (a-b), so Switch on Failure is superior to the single observation strategy if a + b > 1  They yield the same success probability if a + b =1, and Reject on Failure is inferior if a + b < 1. The results of both Stay on Success and Switch on Failure show that, surprisingly, there is no inherent advantage in making two observations rather than one.

**Section III – Discussion and Outlook**

It is shown in ([3]) that it in many situations one need not have the possibility of an equiprobable choice of which experiment to observe in order to successfully guess which experiment will produce the larger numerical outcome.  It is also shown in ([4]) that the Blackwell-Cover Guessing Strategy can be used to separate coins with differing heads probabilities into subsets of unseen coins with different means of heads probabilities.  One can anticipate that these results are capable of greater generalization using the Blackwell-Cover Guessing Strategy outlined in this paper.

**References**

[1] D. Blackwell (1951), *On the translation parameter problem for discrete variables*, Annals of Mathematical Statistics 22: 391–399. doi10.1214/aoms/1177729585; MR0043418

[2] T. Cover (1987), *Pick the largest number*, In: Open Problems in Communication and Computation, T. Cover and B. Gopinath eds, Springer, p. 152. doi10. 10071/ 978-1-4612-4808-8

[3] Stein J. (2026) Using the Cover Guessing Game Strategy to Predict the Results of Stochastic Processes **arXiv:2608.01300**

[4] Stein J. [2026] Using the Blackwell-Cover Guessing Strategy to Separate Unobserved Bernoulli Trials – II **arXiv:2610.36016v1**